\UseRawInputEncoding 
\documentclass[11pt,a4paper]{article}
\usepackage{amsfonts}
\usepackage{amssymb}
\usepackage{mathrsfs}
\usepackage{amsmath}
\usepackage{booktabs}
\usepackage{epsf,epsfig,amsfonts,amsgen,indentfirst}
\usepackage{amsmath,amstext,amsbsy,amsopn,amsthm,bbding,wasysym}
\usepackage{multicol,mathdots}
\usepackage{subfigure}
\usepackage[numbers,sort&compress]{natbib}
\allowdisplaybreaks
\DeclareMathOperator{\Span}{span}

\usepackage{amssymb}

\newcommand{\cA}{\mathcal{A}}

\newtheorem{theorem}{Theorem}[section]
\newtheorem{corollary}{Corollary}[section]

\newtheorem{lemma}[theorem]{Lemma}

\newtheorem{definition}{Definition}[section]

\newtheorem{problem}[theorem]{Problem}

\begin{document}

\title
{\LARGE \textbf{Upper bounds on all-terminal reliability polynomials in terms of cycle basis lengths\thanks{ supported by NSFC (No. 12261071), NSF of Qinghai Province (No. 2025-ZJ-902T).} }}
\author{Tingzeng Wu$^{a,b}$\thanks{{Corresponding author.\newline
\emph{E-mail addresses}: mathtzwu@163.com, 19571327830@163.com, jianxuanluo@163.com, hjlai2015@hotmail.com
}}, Xueqing Li$^{a}$, Jianxuan Luo$^{c}$,  Hong-Jian Lai$^{d}$\\
{\small $^{a}$ School of Mathematics and Statistics, Qinghai Minzu University, }\\
{\small Xining, Qinghai 810007, P.R.~China} \\
{\small $^{b}$Qinghai Institute of Applied Mathematics, Xining, Qinghai, 810007, P.R.~China}\\
{\small $^{c}$ School of Mathematical Sciences, Xiamen University,}\\
{\small  Xiamen, Fujian 361005,  P.R.~China} \\
{\small $^{d}$ School of Mathematics and System Sciences, Guangdong Polytechnic Normal University,}\\
{\small  Guangzhou, Guangdong 510665, P.R.~China}}

\date{}
\maketitle
\noindent {\bf Abstract:}
Let $G=(V(G),E(G))$ be a connected multigraph with $n$ vertices and $m$ edges. 
The all-terminal reliability polynomial of $G$ can be written as $R_G(q)=(1-q)^{n-1}\sum\limits_{i=0}^{d}H_iq^i$ 
where $d=m-n+1$ is the cycle rank, each edge fails independently with probability $q\in[0,1]$, and $H_0,\ldots,H_d$ are the associated coefficients.
A cycle basis, which consists of cycles that form a basis of the cycle space of a graph, is an important tool for characterizing its cycle structure. The structural properties of cycle bases, together with the associated length optimization problems and construction algorithms, have been extensively studied. In this context, Colbourn posed the following open problem, which concerns upper bounds on the reliability coefficients in terms of cycle-basis lengths: 
For all-terminal reliability, is it true that if $G$ has a cycle basis containing cycles of lengths $c_1,\ldots,\allowbreak c_{m-n+1}$, then $H_i \le\sum\limits_{\substack{I\subseteq\{1,\ldots,m-n+1\},|I|=i}}\prod\limits_{j\in I}(c_j-1)$? In this paper, we establish a result that resolves this open problem.

\smallskip


\section{Introduction}\label{sec:introduction}
Let $G=(V(G),E(G))$ be a connected multigraph with $n$ vertices and $m$ edges. Its all-terminal reliability polynomial $R_G(q)$ gives the probability that the spanning subgraph consisting of the operational edges is connected. We assume that all vertices remain operational and that each edge fails independently with probability $q$.
A cycle is a nonempty closed walk in which no edge is repeated and no vertex is repeated except that the initial and terminal vertices coincide. In particular, a single loop is permitted as a cycle of length $1$, and two parallel edges are permitted as a cycle of length $2$. When computing vertex degrees, each loop is counted twice.
For an edge set $S\subseteq E(G)$, let $G-S$ denote the spanning subgraph obtained by deleting the edges in $S$ while retaining all vertices. An edge set $S$ is called a deletable edge set if $G-S$ is connected. An edge $e$ is called a bridge if $G-\{e\}$ is disconnected, and the number of bridges is denoted by $b$.
For a nonnegative integer $h$, let $[h]=\{1,\ldots,h\}$, with $[0]=\varnothing$. For a finite set $S$, the notation $|S|$ denotes its cardinality.

The all-terminal reliability polynomial of the graph is defined by
\begin{equation}\label{eq:H}
R_G(q)=(1-q)^{n-1}\sum_{i=0}^{d}H_iq^i,
\end{equation}
where the coefficients $H_i$ are uniquely determined by \eqref{eq:H}.
The sequence $(H_0,\ldots,H_d)$ is called the $H$-sequence of the all-terminal reliability polynomial, with the convention that $H_i=0$ for $i>d$.
Let $G_0$ denote the graph obtained from $G$ by deleting all loops, and let $L$ be the set of loops of $G$, with $\ell=|L|$. Since deleting loops does not affect connectivity, $G_0$ remains connected and $R_G(q)=R_{G_0}(q)$. It follows from \eqref{eq:H} and the uniqueness of the coefficients that the $H$-sequence of $G$ is obtained by appending $\ell$ zeros to the $H$-sequence of $G_0$. In particular, $H_i=0$ whenever $d-\ell<i\le d$.

Research on upper and lower bounds for the all-terminal reliability polynomial has drawn on shellability, network transformations, coefficient inequalities, and extremal graphs. Ball and Provan \cite{bp1982} exploited the shellable structure of the family of deletable edge sets to derive upper and lower bounds on the all-terminal reliability function, incorporating information on the number of spanning trees and the size and number of minimum cuts. Colbourn and Harms \cite{ch1988} employed linear programming to improve the reliability estimates provided by the Ball--Provan and Lomonosov--Polesskii bounds. Brown et al. \cite{bcd1993} investigated the effects of network transformations on the $H$-coefficients of the all-terminal reliability polynomial and established efficiently computable coefficientwise upper and lower bounds, which they incorporated as additional constraints to strengthen the Ball--Provan method. For $H$-sequences arising from graphs, Brown and Colbourn \cite{browncolbourn} established several non-Stanley-type inequalities, which provide new coefficient constraints for reliability estimation.
Brown et al. \cite{bcn2009} exploited the connection between chip-firing processes on graphs and reliability $H$-sequences, together with structural information about the corresponding monomial representations, to improve the all-terminal reliability bounds obtained by the Ball--Provan method. For sparse graphs, Romero \cite{romero2022} considered connected simple graphs with prescribed numbers of vertices and edges whose cycle rank is at most $4$. Using uniformly least reliable graphs and previously known uniformly most reliable graphs, Romero obtained optimal universal upper and lower bounds on reliability for this class of graphs.
Building on earlier results, Colbourn \cite{Colbourn1993} investigated the coefficients of the all-terminal reliability polynomial and posed the following open problem.
\begin{problem}\rmfamily\normalfont\cite{Colbourn1993}\label{conj:reconstruction}
For all-terminal reliability, is it true that if $G$ has a cycle basis containing cycles of lengths $c_1,\ldots,\allowbreak c_{m-n+1}$, then
\[
H_i \le
\sum_{\substack{I\subseteq\{1,\ldots,m-n+1\},|I|=i}}
\prod_{j\in I}(c_j-1)?
\]
\end{problem}

Define a cactus graph to be a connected graph in which each edge belongs to at most one cycle, and let $\mathbb Z_2$ denote the field with two elements.
We investigate this problem and obtain the following theorem:

\begin{theorem}\normalfont\label{thm:main}
Let $\mathcal B=\{C_1,\ldots,C_d\}$ be an arbitrary cycle basis of $G$ over $\mathbb Z_2$, and let $c_j=|E(C_j)|\ge1$ denote the length of $C_j$.
Then, for every $0\le i\le d$, the coefficients in \eqref{eq:H} satisfy
\begin{equation}\label{eq:main}
0\le H_i\le
\sum_{\substack{I\subseteq[d],|I|=i}}
\prod_{j\in I}(c_j-1).
\end{equation}
Where, equality holds in the second inequality simultaneously for all $i$ if and only if $G$ is a cactus graph.\\
Remark. For $i>d-\ell$, equality holds in the first inequality if and only if the connected multigraph obtained from $G$ by deleting all loops has cycle rank at most $i-1$. For $0\le i\le d-\ell$, we have $H_i>0$. In particular, if $G$ is loopless, then the first inequality is strict throughout the stated range.
\end{theorem}

Brown and Colbourn \cite{browncolbourn}, in their study of upper and lower bounds on network reliability, used an expansion of the all-terminal reliability polynomial according to the size of the failed edge set, namely the $F$-form.
Huh \cite{huh} also used this representation in studying the log-concavity of matroid $h$-vectors and its applications to network reliability.
Thus, \eqref{eq:H} can be rewritten as follows:
\begin{equation}\label{eq:R}
R_G(q)=\sum_{i=0}^{d}F_iq^i(1-q)^{m-i}.
\end{equation}
Here, $F_i$ is the number of deletable edge sets containing $i$ edges, and $F_i=0$ for $i>d$.
For example, $F_0=1$, $F_1$ is the number of non-bridge edges in $G$, and $F_d$ is the number of spanning trees of $G$.
A spanning tree is a connected acyclic subgraph that contains all vertices and has one fewer edge than vertices; let $\tau(G)$ denote the number of spanning trees of $G$.
From \eqref{eq:H} and \eqref{eq:R}, we obtain
\begin{equation}\label{eq:basic}
H_0=1,\qquad H_1=n-1-b,\qquad
\sum_{i=0}^{d}H_i=F_d=\tau(G).
\end{equation}

The above discussion and Theorem \ref{thm:main} yield the following result.

\begin{corollary}\normalfont\label{cor:consequences}
Under the assumptions of Theorem \ref{thm:main}, for every $0\le i\le d$, we have
\[
F_i\le
\sum_{\substack{I\subseteq[d],|I|=i}}
\prod_{j\in I}c_j.
\]
Equality always holds when $i=0$; for each $1\le i\le d$, equality holds if and only if $G$ is a cactus graph.

In particular, setting $i=d$ in \eqref{eq:basic} yields
\begin{equation}\label{eq:tau}
\tau(G)=F_d\le\prod_{j=1}^{d}c_j,
\end{equation}
with equality if and only if $G$ is a cactus graph.\\
Remark. Bortoletto et al. \normalfont\cite{blm2025} previously established \eqref{eq:tau} in their study of the geometric structure of periodic timetables.
\end{corollary}

The paper is organized as follows. Section \ref{sec:prelim} introduces cycle spaces, independent sets of matroids, and the associated quotient algebras, and presents the lemmas that are needed for the subsequent proofs. In Section \ref{sec:main}, we establish upper bounds on the coefficients of the all-terminal reliability polynomial in terms of the lengths of the cycles in an arbitrary binary cycle basis, from which we further derive upper bounds on the numbers of deletable edge sets and spanning trees. Section \ref{sec:conclusion} summarizes the main results and poses further related questions.

\section{Preliminaries}\label{sec:prelim}
In this section, we retain the notation introduced in Section 1 and present the definitions and lemmas that are needed to prove the main theorem.

Let $\Span_{\mathbb Z_2}$ and $\dim_{\mathbb Z_2}$ denote linear span and vector space dimension over the field $\mathbb Z_2$, respectively. Let $\mathbb Z_2^{E(G)}$ denote the vector space over $\mathbb Z_2$ whose coordinates are indexed by the edge set $E(G)$. For any $S\subseteq E(G)$, let $\chi_S\in\mathbb Z_2^{E(G)}$ denote the characteristic vector of $S$, whose coordinate at an edge $e\in E(G)$ is given by
$$
(\chi_S)_e=
\begin{cases}
1,& e\in S,\\
0,& e\notin S.
\end{cases}
$$
Fix a spanning tree $T$ of $G$. For each non-tree edge $e\in E(G)\setminus E(T)$, let $C_T^e$ denote the fundamental cycle associated with $e$ with respect to $T$. If $e$ is not a loop, then $C_T^e$ consists of $e$ together with the unique path in $T$ joining its endpoints; if $e$ is a loop, then $C_T^e$ is the cycle of length $1$ consisting solely of $e$.

We next introduce the binary cycle space and its basic properties, and describe a matching between the cycles in a cycle basis and the non-tree edges of a spanning tree.
\begin{definition}\normalfont{\cite{kavitha2009}}
The binary cycle space of $G$ is defined by
$$
\mathcal C_{\mathbb Z_2}(G)
=
\Span_{\mathbb Z_2}
\{\chi_C:C\text{ is a cycle of }G\},
$$
where $\chi_C$ denotes the characteristic vector of the edge set of $C$.

A family of cycles $\mathcal B=\{C_1,\ldots,C_d\}$ is called a binary cycle basis of $G$ if the corresponding vectors $\{\chi_{C_1},\ldots,\chi_{C_d}\}$ form a basis of $\mathcal C_{\mathbb Z_2}(G)$.
\end{definition}

Since the only cycle that contains a given loop is the cycle consisting of that loop alone, every binary cycle basis must contain the cycle of length $1$ formed by each loop. Removing these cycles leaves a binary cycle basis of the loopless graph $G_0$.

\begin{lemma}\normalfont{\cite{kavitha2009}}
\label{lem:prelim-fundamental-basis}
The dimension of the binary cycle space of a connected graph $G$ equals its cycle rank, that is,
$$
\dim_{\mathbb Z_2}\mathcal C_{\mathbb Z_2}(G)=d=m-n+1.
$$
\end{lemma}

Moreover, for any spanning tree $T$ of $G$, the characteristic vectors of the fundamental cycles with respect to $T$ form a basis of the binary cycle space.

\begin{lemma}\normalfont{\cite{vanderMeulen2015}}
\label{lem:cycle-basis-matching}
Let $S=\{e_1,\ldots,e_d\}\subseteq E(G)$ be such that $G-S$ is a spanning tree and let $\mathcal B=\{C_1,\ldots,C_d\}$ be a binary cycle basis of $G$. 
Then a bijection
$$
f:S\longrightarrow\mathcal B,
\qquad e\in E(f(e))\quad(e\in S).
$$
\end{lemma}

To establish the connection between the deletable edge sets of a graph and its cographic matroid, we introduce the basic notions of matroids, including independent sets, bases, and rank.
\begin{definition}\normalfont{\cite{oxley2011}} \label{mathcal}
A matroid $M$ is an ordered pair $(E,\mathcal{I})$ consisting of a finite set $E$ and a
collection $\mathcal{I}$ of subsets of $E$ having the following three properties:
\begin{enumerate}\normalfont
\item[(i)] $\varnothing\in\mathcal{I}$;
\item[(ii)] If $I\in\mathcal{I}$ and $I'\subseteq I$, then $I'\in\mathcal{I}$;
\item[(iii)] If $I_1$ and $I_2$ are in $\mathcal{I}$ and $|I_1|<|I_2|$, then there is an element $e$ of $I_2-I_1$ such that $I_1\cup\{e\}\in\mathcal{I}$.
\end{enumerate}

We shall call (i) and (ii) the hereditary and independence augmentation properties.
The members of $\mathcal{I}$ are the independent sets of $M$ , and $E$ is the ground set of $M$ . We shall often write $\mathcal{I}(M)$ for $\mathcal{I}$ and $E(M)$ for $E$, particularly. A subset of $E$ that is not in $E$ is called dependent.
\end{definition}

We call a maximal independent set in $M$ a basis or a base of $M$. 

\begin{lemma}\normalfont{\cite{oxley2011}}
If $B_1$ and $B_2$ are bases of a matroid $M$, then $|B_1|=|B_2|$.
\end{lemma}
All bases of a matroid $M$ have the same size; this common size is called the rank of $M$ and is denoted by $r(M)$.

Let $W$ be a finite set, and let $\Gamma$ be a nonempty family of subsets of $W$. We call $\Gamma$ a simplicial complex on $W$ if $F'\in\Gamma$ whenever $F\in\Gamma$ and $F'\subseteq F$. Each member of $\Gamma$ is called a face.
By axioms~(i) and~(ii) of Definition~\ref{mathcal}, $\mathcal I(M)$ contains $\varnothing$ and is closed under taking subsets. Thus, $\mathcal I(M)$ is a simplicial complex, which is called the independence complex of the matroid $M$.

To study the coefficients $H_i$ using algebraic methods, we introduce a graded quotient algebra associated with $G$.

\begin{definition}\normalfont{\cite{bcw1996}}
Let $P=\mathbb Z_2[x_e:e\in E(G)]$, where each variable $x_e$ has degree $1$. For any $X\subseteq E(G)$, write $x^X=\prod\limits_{e\in X}x_e$, with $x^\varnothing=1$. A monomial is called squarefree if the exponent of each variable is either $0$ or $1$. Define
$$
\begin{aligned}
I_{\mathrm{cut}}
&=
\left\langle
x^S:S\subseteq E(G),\ G-S\text{ is disconnected}
\right\rangle,\\
J_{\mathrm{cyc}}
&=
\left\langle
\sum_{e\in E(C)}x_e:
C\text{ is a cycle of }G
\right\rangle.
\end{aligned}
$$
Both $I_{\mathrm{cut}}$ and $J_{\mathrm{cyc}}$ are homogeneous ideals.
Define the graded quotient algebra
\begin{equation}\label{eq:algebra}
\cA(G)
=
P/(I_{\mathrm{cut}}+J_{\mathrm{cyc}})
=
\bigoplus_{i\ge0}\cA(G)_i,
\end{equation}
where $\cA(G)_i$ is the homogeneous component of degree $i$. Let $\bar x_e$ denote the image of $x_e$ in the quotient algebra, and write $\bar x^S=\prod\limits_{e\in S}\bar x_e$, with $\bar x^\varnothing=1$.
\end{definition}
Fix a binary cycle basis $\mathcal B=\{C_1,\ldots,C_d\}$ of $G$. Since $\chi_{C_1},\ldots,\chi_{C_d}$ form a basis of the binary cycle space of $G$, we have
\begin{equation}\label{eq:parameters}
\begin{aligned}
J_{\mathrm{cyc}}
&=(\theta_1,\ldots,\theta_d),\qquad
\theta_j
&=\sum_{e\in E(C_j)}x_e,
\qquad 1\le j\le d.
\end{aligned}
\end{equation}
In particular, $\sum\limits_{e\in E(C_j)}\bar x_e=0$. If $e$ is a bridge of $G$, then $G-\{e\}$ is disconnected, so $x_e\in I_{\mathrm{cut}}$ and hence $\bar x_e=0$.

Let $M(G)$ denote the graphic matroid of $G$, whose independent sets are the edge sets that contain no cycles. Since $G$ is connected, the bases of $M(G)$ are precisely the edge sets of the spanning trees of $G$. Its dual matroid $M(G)^*$ is called the cographic matroid of $G$, whose bases are the complements of the edge sets of the spanning trees of $G$. Consequently, $r(M(G)^*)=m-(n-1)=d$.

By axiom~(ii) of Definition~\ref{mathcal}, every subset of a basis is independent. Moreover, since the ground set is finite, every independent set can be extended to a basis. Thus, an edge set is independent in $M(G)^*$ if and only if it is contained in a basis. For any $S\subseteq E(G)$, we have
$$
\begin{aligned}
S\in\mathcal I(M(G)^*)
&\iff
\text{there exists a spanning tree }T\text{ such that }
S\subseteq E(G)\setminus E(T)\\
&\iff
\text{there exists a spanning tree }T\text{ such that }
E(T)\subseteq E(G)\setminus S\\
&\iff
G-S\text{ is connected}.
\end{aligned}
$$
Therefore, the independent sets of $M(G)^*$ are precisely the deletable edge sets of $G$, and its independence complex is $\mathcal I(M(G)^*)=\{S\subseteq E(G):G-S\text{ is connected}\}$.
The graded quotient ring $\mathbb Z_2[M(G)^*]=P/\langle x^S:S\subseteq E(G),\ S\notin\mathcal I(M(G)^*)\rangle$ is called the Stanley--Reisner ring of the independence complex of the cographic matroid $M(G)^*$.
Define $h_0,\ldots,h_d$ by
\begin{equation}\label{eq:sumIMG}
\sum\limits_{S\in\mathcal I(M(G)^*)}t^{|S|}(1-t)^{d-|S|}=\sum\limits_{i=0}^{d}h_it^i.
\end{equation}
The coefficient vector $(h_0,\ldots,h_d)$ is called the $h$-vector of $\mathcal I(M(G)^*)$.

\begin{lemma}\normalfont{\cite{bcw1996}}
\label{lem:bcw-dimension}
Let $T$ be a spanning tree of $G$, and set
$$
\Theta=
\left\{
\sum_{f\in E(C_T^e)}x_f:e\in E(G)\setminus E(T)
\right\}.
$$
Then the quotient algebra admits the graded decomposition
$$
\mathbb Z_2[M(G)^{*}]/\langle\Theta\rangle
=
R'_0\oplus R'_1\oplus\cdots\oplus R'_d,
$$
where $\langle\Theta\rangle$ denotes the ideal generated by the images of the linear forms in $\Theta$ in $\mathbb Z_2[M(G)^*]$, and $R'_i$ denotes the homogeneous component of degree $i$ of $\mathbb Z_2[M(G)^*]/\langle\Theta\rangle$, namely, the vector space consisting of the images of all homogeneous polynomials of degree $i$ in this quotient algebra. The dimensions of these homogeneous components form the $h$-vector of $M(G)^*$; that is, $\dim_{\mathbb Z_2}R'_i=h_i$.
\end{lemma}

Since the bases of $M(G)^*$ are the complements of the edge sets of spanning trees, its independent sets are precisely the deletable edge sets. Hence, for each $0\le k\le d$, there are exactly $F_k$ independent sets of cardinality $k$. By \eqref{eq:H}, \eqref{eq:R}, and the identity $m=n-1+d$, we obtain
$$
\sum_{i=0}^{d}H_iq^i
=
\sum_{k=0}^{d}F_kq^k(1-q)^{d-k}.
$$
Equation \eqref{eq:sumIMG} also gives
$$
\sum_{i=0}^{d}h_iq^i
=
\sum_{k=0}^{d}F_kq^k(1-q)^{d-k}.
$$
Therefore, $\sum\limits_{i=0}^{d}H_iq^i=\sum\limits_{i=0}^{d}h_iq^i$. Comparing coefficients yields $H_i=h_i$ for every $0\le i\le d$. Thus, $(H_0,\ldots,H_d)$ is precisely the $h$-vector of the independence complex of $M(G)^*$.

\begin{lemma}\normalfont{\cite{bcw1996}}
\label{lem:bcw-squarefree}
The images of the squarefree monomials $\{x^S:S\subseteq E(G)\}$ in the quotient algebra $P/(I_{\mathrm{cut}}+J_{\mathrm{cyc}})$ span this algebra.
\end{lemma}
When $G$ is loopless, the assertion follows directly from \cite{bcw1996}. If $G$ contains loops, then each loop $e$ forms a cycle by itself, so $x_e\in J_{\mathrm{cyc}}$ and hence $\bar x_e=0$. Since deleting loops does not affect the connectivity of $G$, we have $\cA(G)\cong\cA(G_0)$. Applying Lemma~\ref{lem:bcw-squarefree} to $G_0$ therefore yields the desired spanning result for $\cA(G)$.

To identify the indices for which $H_i$ is positive or zero, we state the following result on the $h$-vector of a matroid independence complex.
\begin{lemma}\normalfont{\cite{browncolbourn}}
\label{lem:h-zero}
Let $M$ be a finite matroid of rank $r$ whose independence complex $\mathcal I(M)$ has $h$-vector $(h_0,\ldots,h_r)$. If $M$ has exactly $a$ elements that belong to every basis, then
$$
h_i>0\quad(0\le i\le r-a),\qquad
h_i=0\quad(r-a<i\le r).
$$
\end{lemma}

We next recall the third isomorphism theorem for rings, which will be used to handle successive quotients.
\begin{lemma}\normalfont{\cite{noether1927}}
\label{lem:third-isomorphism}
Let $R$ be an associative ring, and let $I,J$ be two-sided ideals of $R$ with $I\subseteq J$. Then $J/I$ is a two-sided ideal of the quotient ring $R/I$, and the map
$$
\begin{aligned}
\varphi:(R/I)/(J/I)&\longrightarrow R/J,\\
(r+I)+(J/I)&\longmapsto r+J
\end{aligned}
$$
is a ring isomorphism. Consequently, $(R/I)/(J/I)\cong R/J$.
\end{lemma}

We next relate the quotient algebra defined above to the coefficients $H_i$ and identify these coefficients with the dimensions of its homogeneous components.

\begin{lemma}\label{prop:known}
For every $i\ge0$, we have
$$
H_i=
\begin{cases}
\dim_{\mathbb Z_2}\cA(G)_i,&0\le i\le d,\\
0,&i>d,
\end{cases}
$$
and
$$
\cA(G)_i=
\begin{cases}
\Span_{\mathbb Z_2}
\left\{
\bar x^S:S\subseteq E(G),\ |S|=i,\ G-S\text{ is connected}
\right\},&0\le i\le d,\\
\{0\},&i>d.
\end{cases}
$$
\end{lemma}

\begin{proof}
First, since the independent sets of $M(G)^*$ are precisely the deletable edge sets of $G$, for any $S\subseteq E(G)$ we have $S\notin\mathcal I(M(G)^*)$ if and only if $G-S$ is disconnected. Thus, the Stanley--Reisner ideal of the independence complex of $M(G)^*$ is
$$
\begin{aligned}
\left\langle
x^S:S\subseteq E(G),\ S\notin\mathcal I(M(G)^*)
\right\rangle
=
\left\langle
x^S:S\subseteq E(G),\ G-S\text{ is disconnected}
\right\rangle
=I_{\mathrm{cut}}.
\end{aligned}
$$
Consequently, $\mathbb Z_2[M(G)^*]=P/I_{\mathrm{cut}}$.

Next, fix a spanning tree $T$ of $G$. Let $\theta_e=\sum\limits_{f\in E(C_T^e)}x_f$, so that $\Theta=\{\theta_e:e\in E(G)\setminus E(T)\}$. Let $\langle\Theta\rangle_P$ denote the ideal generated by the linear forms in $\Theta$ in the polynomial ring $P$. We now show that $\langle\Theta\rangle_P=J_{\mathrm{cyc}}$.

We first prove that $J_{\mathrm{cyc}}\subseteq\langle\Theta\rangle_P$. Let $C$ be an arbitrary cycle of $G$. Since the characteristic vectors of the fundamental cycles $C_T^e$ form a basis of the binary cycle space, there exist coefficients $a_e\in\mathbb Z_2$ such that $\chi_C=\sum\limits_{e\in E(G)\setminus E(T)}a_e\chi_{C_T^e}$. For any $f\in E(G)$, comparing the coordinates indexed by $f$ on both sides yields
\begin{equation}\label{eq:chi}
(\chi_C)_f
=
\sum_{e\in E(G)\setminus E(T)}
a_e(\chi_{C_T^e})_f.
\end{equation}
Since $(\chi_C)_f=1$ when $f\in E(C)$ and $(\chi_C)_f=0$ otherwise, we have
\begin{equation}\label{eq:sum}
\sum_{f\in E(G)}(\chi_C)_f x_f=\sum_{f\in E(C)}1\cdot x_f+\sum_{f\in E(G)\setminus E(C)}0\cdot x_f=\sum_{f\in E(C)}x_f.
\end{equation}
Combining \eqref{eq:chi} and \eqref{eq:sum}, we obtain
$
\sum\limits_{f\in E(C)}x_f
=\sum\limits_{e\in E(G)\setminus E(T)}a_e\theta_e.
$
Each $\theta_e$ belongs to $\langle\Theta\rangle_P$, which is closed under addition and multiplication by elements of $P$. Hence, the above linear combination also belongs to $\langle\Theta\rangle_P$, and therefore $\sum\limits_{f\in E(C)}x_f\in\langle\Theta\rangle_P$. Since $C$ is arbitrary, $\langle\Theta\rangle_P$ contains the linear form associated with every cycle of $G$. As $J_{\mathrm{cyc}}$ is precisely the ideal generated by all these linear forms, we conclude that $J_{\mathrm{cyc}}\subseteq\langle\Theta\rangle_P$.

Conversely, every fundamental cycle $C_T^e$ is itself a cycle of $G$. Since $J_{\mathrm{cyc}}$ is the ideal generated by the linear forms associated with all cycles of $G$, we have $\theta_e=\sum\limits_{f\in E(C_T^e)}x_f\in J_{\mathrm{cyc}}$ for every non-tree edge $e\in E(G)\setminus E(T)$. Thus, all linear forms in $\Theta$ belong to $J_{\mathrm{cyc}}$. Since $J_{\mathrm{cyc}}$ is an ideal, it also contains all polynomial multiples of these linear forms and all finite sums of such multiples. By the definition of $\langle\Theta\rangle_P$ as the ideal generated by the linear forms in $\Theta$ in the polynomial ring $P$, it follows that $\langle\Theta\rangle_P\subseteq J_{\mathrm{cyc}}$.

Combining the two inclusions yields $\langle\Theta\rangle_P=J_{\mathrm{cyc}}$.

Next, let $\langle\Theta\rangle$ denote the ideal in $P/I_{\mathrm{cut}}$ generated by the images of the linear forms in $\Theta$. Since every element of $I_{\mathrm{cut}}$ maps to zero in the quotient ring, adding an element of $I_{\mathrm{cut}}$ to an element of $\langle\Theta\rangle_P$ does not change its image. Hence, $\langle\Theta\rangle=\bigl(I_{\mathrm{cut}}+\langle\Theta\rangle_P\bigr)/I_{\mathrm{cut}}$. Therefore, by Lemma~\ref{lem:third-isomorphism} and \eqref{eq:algebra}, we have
$$
\begin{aligned}
\frac{\mathbb Z_2[M(G)^*]}{\langle\Theta\rangle}
=
\frac{\mathbb Z_2[M(G)^*]}{\bigl(I_{\mathrm{cut}}+\langle\Theta\rangle_P\bigr)/I_{\mathrm{cut}}}
=
\frac{P/I_{\mathrm{cut}}}{(I_{\mathrm{cut}}+\langle\Theta\rangle_P)/I_{\mathrm{cut}}}
\cong
\frac{P}{I_{\mathrm{cut}}+\langle\Theta\rangle_P}
=
\frac{P}{I_{\mathrm{cut}}+J_{\mathrm{cyc}}}
=
\cA(G).
\end{aligned}
$$
By Lemma~\ref{lem:bcw-dimension}, it follows that
$$
\dim_{\mathbb Z_2}\cA(G)_i
=
\dim_{\mathbb Z_2}R'_i=h_i=H_i.
$$
For $i>d$, we have $\cA(G)_i=\{0\}$ and $H_i=0$, so $\dim_{\mathbb Z_2}\cA(G)_i=H_i=0$.

Finally, we prove the spanning assertion for each homogeneous component. For $i\geq0$, let $u\in\cA(G)_i$ be arbitrary. By Lemma~\ref{lem:bcw-squarefree}, there exist coefficients $a_S\in\mathbb Z_2$ such that $u=\sum\limits_{S\subseteq E(G)}a_S\bar x^S$. Since $I_{\mathrm{cut}}+J_{\mathrm{cyc}}$ is a homogeneous ideal, $\cA(G)$ is a graded algebra in which each $\bar x^S$ is homogeneous of degree $|S|$. Taking the homogeneous component of degree $i$ on both sides yields $u=\sum\limits_{\substack{S\subseteq E(G),|S|=i}}a_S\bar x^S$. If $G-S$ is disconnected, then $x^S\in I_{\mathrm{cut}}$ and hence $\bar x^S=0$. Omitting these zero terms, we obtain
$$
u=
\sum_{\substack{
S\subseteq E(G),\ |S|=i,\\G-S\text{ is connected}
}}
a_S\bar x^S.
$$
Therefore,
$\cA(G)_i\subseteq\Span_{\mathbb Z_2}\left\{\bar x^S:S\subseteq E(G),\ |S|=i,\ G-S\text{ is connected}\right\}$.
Conversely, each $\bar x^S$ in this set belongs to $\cA(G)_i$, so their span is also contained in $\cA(G)_i$; that is,
$\cA(G)_i\supseteq\Span_{\mathbb Z_2}\left\{\bar x^S:S\subseteq E(G),\ |S|=i,\ G-S\text{ is connected}\right\}$.
Thus,
$$
\cA(G)_i=
\Span_{\mathbb Z_2}
\left\{
\bar x^S:S\subseteq E(G),\ |S|=i,\ G-S\text{ is connected}
\right\}.
$$

If $G-S$ is connected, then it has at least $n-1$ edges, so $m-|S|\ge n-1$ and hence $|S|\le m-n+1=d$. Thus, when $i>d$, there is no deletable edge set of cardinality $i$, and the set of monomials in the spanning expression above is empty. Since the span of the empty set is the zero space, we obtain $\cA(G)_i=\Span_{\mathbb Z_2}(\varnothing)=\{0\}$.
\end{proof}

To estimate the dimensions of higher-degree homogeneous components using dimension bounds for the degree-one subspaces associated with the cycles in a cycle basis, we state two standard dimension estimates for products and sums of subspaces. Inequality \eqref{eq:product-dimension} follows directly from the tensor product basis theorem in \cite{lang2002}, while \eqref{eq:sum-dimension} follows from the dimension estimate for sums of subspaces in \cite{bourbaki1989}.

\begin{lemma}\normalfont\label{lem:product}
Let $A$ be a nonzero unital commutative graded algebra over a field $K$, let $s$ be a nonnegative integer, and let $V_1,\ldots,V_s\subseteq A_1$ be finite-dimensional subspaces. Denote by $V_1\cdots V_s$ the $K$-vector space spanned by all products $v_1\cdots v_s$ with $v_j\in V_j$. Then
\begin{equation}\label{eq:product-dimension}
\dim_K(V_1\cdots V_s)
\le \prod_{j=1}^{s}\dim_K V_j.
\end{equation}
If $\{W_\alpha\}_{\alpha\in\Lambda}$ is a finite family of finite-dimensional subspaces of a common $K$-vector space, then
\begin{equation}\label{eq:sum-dimension}
\dim_K\left(\sum_{\alpha\in\Lambda}W_\alpha\right)
\le \sum_{\alpha\in\Lambda}\dim_K W_\alpha.
\end{equation}
By convention, the empty product space for $s=0$ is $K\cdot1_A$, and the sum of an empty family of subspaces is $\{0\}$. The empty numerical product and sum are taken to be $1$ and $0$, respectively.
\end{lemma}

\section{Proof of Theorem \ref{thm:main}}\label{sec:main}
\begin{proof}[\rm\textbf{Proof of Theorem \ref{thm:main}}]
By Lemma~\ref{prop:known}, for every $0\le i\le d$, we have $H_i=\dim_{\mathbb Z_2}\cA(G)_i\ge0$. Moreover, $\cA(G)_i$ is spanned by the squarefree monomials $\bar x^S$ for which $|S|=i$ and $G-S$ is connected. Thus, to prove \eqref{eq:main}, it suffices to show that
$$
\dim_{\mathbb Z_2}\cA(G)_i
\le
\sum_{\substack{I\subseteq[d],|I|=i}}
\prod_{j\in I}(c_j-1).
$$
To this end, we first match the edges of each deletable edge set to distinct cycles in the given cycle basis and then show that the corresponding monomial lies in the product space determined by these cycles.

We first show that, for every deletable edge set $S\subseteq E(G)$, there exists an injection

$$
\varphi:S\longrightarrow[d],
\qquad e\in E(C_{\varphi(e)})\quad(e\in S).
$$

If $S=\varnothing$, the empty map suffices. Henceforth, assume that $S\ne\varnothing$. Since $G-S$ is connected, choose a spanning tree $T$ of $G-S$. Because all vertices are retained when edges are deleted, $T$ is also a spanning tree of $G$. Moreover, $E(T)\subseteq E(G)\setminus S$, so $S\subseteq E(G)\setminus E(T)$. By Lemma~\ref{lem:cycle-basis-matching}, there exists a bijection
$$
f:E(G)\setminus E(T)\longrightarrow\mathcal B,
\qquad e\in E(f(e)).
$$
For each $e\in S$, define $\varphi(e)=j$ whenever $f(e)=C_j$. Since $\mathcal B=\{C_1,\ldots,C_d\}$, each $e\in S$ determines a unique index $\varphi(e)\in[d]$. The injectivity of $f$ implies that $\varphi$ is also injective. For every $e\in S$, we have $e\in E(f(e))=E(C_{\varphi(e)})$.

Next, for each $j\in[d]$, define
$$
V_j=\Span_{\mathbb Z_2}\{\bar x_e:e\in E(C_j)\}
\subseteq\cA(G)_1.
$$
The space $V_j$ is spanned by the images of the $c_j$ edge variables corresponding to the edges of $C_j$. By \eqref{eq:parameters}, we have $\sum\limits_{e\in E(C_j)}x_e\in J_{\mathrm{cyc}}$, so $\sum\limits_{e\in E(C_j)}\bar x_e=0$ in the quotient algebra. Since $C_j$ is nonempty, choose any edge $e_0\in E(C_j)$. Over the field $\mathbb Z_2$, we have $\bar x_{e_0}=\sum\limits_{e\in E(C_j)\setminus\{e_0\}}\bar x_e$, which expresses $\bar x_{e_0}$ as a linear combination of the images of the remaining $c_j-1$ edge variables. Thus, removing $\bar x_{e_0}$ from the spanning family leaves a family that still spans $V_j$; that is,
$$
V_j=\Span_{\mathbb Z_2}
\{\bar x_e:e\in E(C_j)\setminus\{e_0\}\}.
$$
Since every vector in $V_j$ can be expressed as a linear combination of these $c_j-1$ vectors, the dimension of $V_j$ is at most $c_j-1$; that is,
\begin{equation}\label{eq:dimU}
\dim_{\mathbb Z_2}V_j\le c_j-1.
\end{equation}
In particular, when $c_j=1$, the cycle $C_j$ consists of a single loop, so $V_j=\{0\}$ and \eqref{eq:dimU} remains valid.

Next, for each $I\subseteq[d]$, let
$$
W_I=\prod_{j\in I}V_j
=\Span_{\mathbb Z_2}
\left\{\prod_{j\in I}v_j:v_j\in V_j\ (j\in I)\right\}.
$$
By Lemma~\ref{lem:product} and \eqref{eq:dimU}, we obtain
\begin{equation}\label{eq:dimW}
\dim_{\mathbb Z_2}W_I
\le\prod_{j\in I}\dim_{\mathbb Z_2}V_j
\le\prod_{j\in I}(c_j-1).
\end{equation}

We now prove that $\cA(G)_i=\sum\limits_{\substack{I\subseteq[d],|I|=i}}W_I$.

We first show that $\cA(G)_i\subseteq\sum\limits_{\substack{I\subseteq[d],|I|=i}}W_I$. Let $S$ be a nonempty deletable edge set of cardinality $i$. There exists an injection $\varphi:S\longrightarrow[d]$ such that $e\in E(C_{\varphi(e)})$ for every $e\in S$. Set $I=\varphi(S)$, so that $\varphi$ induces a bijection between $S$ and $I$. For each $j\in I$, let $e\in S$ be the unique edge for which $\varphi(e)=j$, and set $v_j=\bar x_e$. Since $e\in E(C_j)$, we have $v_j\in V_j$. Thus, the product $\prod\limits_{j\in I}v_j$ is formed by choosing one element from each $V_j$. By the definition $W_I=\prod\limits_{j\in I}V_j$, we obtain $\bar x^S=\prod\limits_{e\in S}\bar x_e=\prod\limits_{j\in I}v_j\in W_I$.
When $S=\varnothing$, we take $I=\varnothing$, and then $\bar x^S=1\in W_\varnothing$, so the assertion remains valid. By Lemma~\ref{prop:known}, $\cA(G)_i$ is spanned by all monomials $\bar x^S$ for which $|S|=i$ and $G-S$ is connected. Each such monomial belongs to some $W_I$ with $|I|=i$. Consequently, every element of $\cA(G)_i$ belongs to the sum of these subspaces, which gives $\cA(G)_i\subseteq\sum\limits_{\substack{I\subseteq[d],|I|=i}}W_I$.

We next prove that $\sum\limits_{\substack{I\subseteq[d],|I|=i}}W_I\subseteq\cA(G)_i$. Since $\cA(G)$ is a graded quotient algebra, its multiplication satisfies $\cA(G)_{r_1}\,\cA(G)_{r_2}\subseteq\cA(G)_{r_1+r_2}$ for all $r_1,r_2\ge0$.
Write $I=\{j_1,\ldots,j_s\}$, where $s=|I|=i\ge1$. For each $l\in\{1,\ldots,s\}$, choose an arbitrary $v_{j_l}\in V_{j_l}$. Since $V_{j_l}\subseteq\cA(G)_1$, repeated application of the multiplication property yields
$$
v_{j_1}\cdots v_{j_s}
\in
\cA(G)_{\underbrace{1+\cdots+1}_{s\text{ terms}}}
=
\cA(G)_s.
$$
Every element of $W_I$ is a linear combination of products of the form $\prod\limits_{j\in I}v_j$, so $W_I\subseteq\cA(G)_s=\cA(G)_i$. When $I=\varnothing$, we have $W_\varnothing=\mathbb Z_2\cdot1\subseteq\cA(G)_0$, so the assertion remains valid. Since each $W_I$ with $|I|=i$ is contained in $\cA(G)_i$, their sum is also contained in $\cA(G)_i$. Hence, $\sum\limits_{\substack{I\subseteq[d],|I|=i}}W_I\subseteq\cA(G)_i$.

Combining the two inclusions, we obtain
\begin{equation}\label{eq:A(G)W}
\cA(G)_i
=
\sum_{\substack{I\subseteq[d],|I|=i}}W_I.
\end{equation}

Finally, applying Lemmas~\ref{prop:known} and~\ref{lem:product}, together with \eqref{eq:dimW} and \eqref{eq:A(G)W}, yields
$$
\begin{aligned}
H_i
=\dim_{\mathbb Z_2}\cA(G)_i
=\dim_{\mathbb Z_2}
\left(\sum_{\substack{I\subseteq[d],|I|=i}}W_I\right)
\le\sum_{\substack{I\subseteq[d],|I|=i}}
\dim_{\mathbb Z_2}W_I
\le\sum_{\substack{I\subseteq[d],|I|=i}}
\prod_{j\in I}(c_j-1).
\end{aligned}
$$
Since $H_i$ is nonnegative, we conclude that
$$
0\le H_i\le
\sum_{\substack{I\subseteq[d],|I|=i}}
\prod_{j\in I}(c_j-1),
$$
where $d=m-n+1$.

We now show that equality holds in the second inequality in \eqref{eq:main} simultaneously for all $0\le i\le d$ if and only if $G$ is a cactus graph.
We first prove necessity. Suppose that $H_i=\sum\limits_{\substack{I\subseteq[d],|I|=i}}\prod_{j\in I}(c_j-1)$ for every $0\le i\le d$. If $d=0$, then $G$ is a tree and hence a cactus graph. Henceforth, assume that $d\ge1$. Taking $i=1$ gives $H_1=\sum\limits_{j=1}^{d}(c_j-1)$. Using $H_1=n-1-b$ and $d=m-n+1$, we obtain
\begin{equation}\label{eq:cj}
\sum_{j=1}^{d}c_j
=H_1+d
=(n-1-b)+(m-n+1)
=m-b.
\end{equation}

First, we show that the union of the edge sets of the basis cycles is precisely the set of all non-bridge edges of $G$. Let $e$ be an arbitrary non-bridge edge of $G$. If $e$ is a loop, then it forms a cycle $C$ by itself. Otherwise, since $G-\{e\}$ is connected, there exists a simple path between the endpoints of $e$ that avoids $e$. This path, together with $e$, forms a cycle $C$.
Since the characteristic vectors of the cycles in $\mathcal B=\{C_1,\ldots,C_d\}$ span the binary cycle space, there exist $a_1,\ldots,a_d\in\mathbb Z_2$ such that $\chi_C=\sum\limits_{j=1}^{d}a_j\chi_{C_j}$. Suppose that $e$ belongs to none of the basis cycles. Then $(\chi_{C_j})_e=0$ for every $j\in[d]$, and hence $(\chi_C)_e=\sum\limits_{j=1}^{d}a_j(\chi_{C_j})_e=0$. This contradicts $(\chi_C)_e=1$, which holds because $e\in E(C)$. Thus, $e$ belongs to at least one basis cycle. Since $e$ was arbitrary, every non-bridge edge belongs to at least one basis cycle.

Choose any edge of a fundamental cycle. If it is a self-loop, deleting it does not affect connectivity; otherwise, after its deletion, its endpoints remain connected via the remaining edges of the cycle. Thus, deleting any edge of a basis cycle leaves $G$ connected, so every edge of a basis cycle is a non-bridge edge.

Thus, the union of the edge sets of the basis cycles is precisely the set of all non-bridge edges of $G$, so $\left|\bigcup\limits_{j=1}^{d}E(C_j)\right|=m-b$. Combining $c_j=|E(C_j)|$ with \eqref{eq:cj}, we obtain
$$
\sum_{j=1}^{d}|E(C_j)|
=
\sum_{j=1}^{d}c_j
=
m-b
=
\left|\bigcup_{j=1}^{d}E(C_j)\right|.
$$
Here, $\sum\limits_{j=1}^{d}|E(C_j)|$ counts the edges separately for each basis cycle, whereas $\left|\bigcup\limits_{j=1}^{d}E(C_j)\right|$ counts each edge in the union exactly once. If two basis cycles share an edge, then that edge is counted at least twice in $\sum\limits_{j=1}^{d}|E(C_j)|$ but only once in $\left|\bigcup\limits_{j=1}^{d}E(C_j)\right|$. Every other edge in the union is counted at least once in the former expression, so $\sum\limits_{j=1}^{d}|E(C_j)|>\left|\bigcup\limits_{j=1}^{d}E(C_j)\right|$, which contradicts the equality above. Therefore, the edge sets of $C_1,\ldots,C_d$ are pairwise disjoint.

Let $C$ be an arbitrary cycle of $G$. Since $\mathcal B=\{C_1,\ldots,C_d\}$ is a binary cycle basis, the edge set of $C$ can be expressed as the symmetric difference of the edge sets of some basis cycles. Since these edge sets are pairwise disjoint, each edge occurs exactly once among the selected cycles. Hence, their symmetric difference coincides with their union.

We next show that only one basis cycle can be selected. Suppose that at least two basis cycles are selected. If the selected cycles are pairwise vertex-disjoint, then the subgraph they form is disconnected, which contradicts the connectedness of $C$. Otherwise, two selected cycles share a vertex, at which each cycle contributes degree $2$, with a loop contributing $2$ to the degree of its vertex. Since these two cycles have disjoint edge sets, their union has degree $4$ at the shared vertex. The subgraph formed by all selected cycles therefore has degree at least $4$ at that vertex, which contradicts the fact that every vertex of $C$ has degree $2$.
Thus, at most one basis cycle can be selected. Since $C$ is nonempty, exactly one is selected, so $C=C_j$ for some $j\in[d]$. As $C$ was arbitrary, the cycles of $G$ are precisely $C_1,\ldots,C_d$. Their edge sets are pairwise disjoint, so every edge of $G$ belongs to at most one cycle. Therefore, $G$ is a cactus graph.

We next prove sufficiency. Suppose that $G$ is a cactus graph. Then the edge sets of all its cycles are pairwise disjoint, so their characteristic vectors are linearly independent. Since these vectors also span the binary cycle space of $G$, they form a basis of this space. As the binary cycle space has dimension $d$, the graph $G$ has exactly $d$ cycles. Thus, the binary cycle basis $\mathcal B=\{C_1,\ldots,C_d\}$ consists of all the cycles of $G$.

Since every non-bridge edge belongs to a cycle and no bridge belongs to any cycle, $E(G)$ is the disjoint union of the set of bridges and the edge sets of $C_1,\ldots,C_d$. Using $c_j=|E(C_j)|$ and $d=m-n+1$, we obtain
$$
m=b+\sum_{j=1}^{d}c_j,
\qquad
n-1=m-d=b+\sum_{j=1}^{d}(c_j-1).
$$

Delete an arbitrary edge from each $C_j$ and retain all bridges. Deleting a loop does not affect connectivity. For each deleted non-loop edge, its endpoints remain connected by the remaining edges of the cycle containing it. Since the cycles have pairwise disjoint edge sets, none of these remaining edges has been deleted, so the resulting spanning subgraph is connected. Moreover, exactly $d$ distinct edges have been deleted, leaving a spanning subgraph with $n$ vertices and $m-d=n-1$ edges. Hence, this subgraph is a spanning tree of $G$.
Conversely, every spanning tree is acyclic and therefore omits at least one edge from each $C_j$. Since the complement of its edge set contains exactly $d$ edges and the $d$ cycles have pairwise disjoint edge sets, exactly one edge is deleted from each cycle, and all bridges are retained. Thus, the spanning trees of $G$ are precisely the spanning subgraphs obtained by deleting one edge from each $C_j$ and retaining all bridges.

It follows that the spanning subgraph consisting of the operational edges is connected if and only if all bridges are operational and at most one edge fails in each $C_j$. Indeed, if the spanning subgraph consisting of the operational edges is connected, then it contains a spanning tree. All failed edges lie in the complement of the edge set of this tree, so no bridge has failed, and each $C_j$ contains at most one failed edge.
Conversely, suppose that all bridges are operational and each $C_j$ contains at most one failed edge. In each cycle containing a failed edge, select that edge; in each cycle containing no failed edge, select an arbitrary edge. Deleting the selected edges yields a spanning tree, all of whose edges are operational. Therefore, the spanning subgraph consisting of the operational edges is connected.

The probability that all bridges are operational is $(1-q)^b$. For each $C_j$, the probability that no edge fails is $(1-q)^{c_j}$, and the probability that exactly one edge fails is $c_jq(1-q)^{c_j-1}$. When $c_j=1$, the condition that at most one edge fails is always satisfied, and its probability is $(1-q)+q=1$, which is consistent with the fact that loops do not affect connectivity. Since edges fail independently, and the set of bridges and the edge sets of the cycles are pairwise disjoint, we have
$$
\begin{aligned}
R_G(q)
&=(1-q)^b
\prod_{j=1}^{d}
\left((1-q)^{c_j}
+c_jq(1-q)^{c_j-1}\right)
&=(1-q)^{n-1}
\prod_{j=1}^{d}\bigl(1+(c_j-1)q\bigr).
\end{aligned}
$$
Comparing this with \eqref{eq:H} gives $\sum\limits_{i=0}^{d}H_iq^i=\prod\limits_{j=1}^{d}\bigl(1+(c_j-1)q\bigr)$. Expanding the product $\prod\limits_{j=1}^{d}\bigl(1+(c_j-1)q\bigr)$ and comparing the coefficients of $q^i$, we obtain
$$
H_i=
\sum_{\substack{I\subseteq[d],|I|=i}}
\prod_{j\in I}(c_j-1),
\qquad 0\le i\le d.
$$
When $d=0$, all empty products are taken to be $1$, so the conclusion remains valid.

We next show that, for $i>d-\ell$, equality holds in the first inequality if and only if the connected multigraph obtained from $G$ by deleting all loops has cycle rank at most $i-1$, and that $H_i>0$ for $0\le i\le d-\ell$. In particular, if $G$ is loopless, then the first inequality is strict throughout the stated range.

Since the bases of $M(G)^*$ are precisely the complements of the edge sets of spanning trees of $G$, an element belongs to every basis of $M(G)^*$ if and only if the corresponding edge belongs to no spanning tree of $G$. No loop belongs to a spanning tree, whereas every non-loop edge is contained in some spanning tree. Thus, $M(G)^*$ has exactly $\ell$ elements that belong to every basis.

Since $r(M(G)^*)=d$ and $(H_0,\ldots,H_d)$ is the $h$-vector of the independence complex of $M(G)^*$, Lemma~\ref{lem:h-zero} yields
$$
H_i>0\quad(0\le i\le d-\ell),
\qquad
H_i=0\quad(d-\ell<i\le d).
$$
The graph obtained by deleting all loops from $G$ remains connected and has cycle rank
$$
(m-\ell)-n+1=d-\ell.
$$
Therefore, for $i>d-\ell$, equality holds in the first inequality if and only if deleting all loops from $G$ leaves a connected multigraph of cycle rank at most $i-1$; for $0\le i\le d-\ell$, we have $H_i>0$. In particular, if $G$ is loopless, then $\ell=0$, so the first inequality is strict for every $0\le i\le d$.
\end{proof}

\begin{proof}[\rm\textbf{Proof of Corollary~\ref{cor:consequences}}]
Set
$
B_k=\sum\limits_{\substack{J\subseteq[d],|J|=k}}\prod\limits_{j\in J}(c_j-1)(0\le k\le d).
$
Substitute $q=t/(1+t)$ into \eqref{eq:H} and \eqref{eq:R}, and multiply both identities by $(1+t)^m$. Using $m=n-1+d$, we obtain
$$
\sum_{i=0}^{d}F_it^i
=
\sum_{k=0}^{d}H_kt^k(1+t)^{d-k}.
$$

Furthermore, we have
$$
\begin{aligned}
\prod_{j=1}^{d}(1+c_jt)
=
\prod_{j=1}^{d}\bigl((1+t)+(c_j-1)t\bigr)
=
\sum_{k=0}^{d}B_kt^k(1+t)^{d-k}.
\end{aligned}
$$
Comparing the coefficients of $t^i$ on both sides yields
$$
\sum_{\substack{I\subseteq[d],|I|=i}}
\prod_{j\in I}c_j-F_i
=
\sum_{k=0}^{i}\binom{d-k}{i-k}(B_K-H_k)
\geq 0.
$$
Therefore, $F_i\le\sum\limits_{\substack{I\subseteq[d],|I|=i}}\prod\limits_{j\in I}c_j$.

Since $B_k\ge H_k$ by Theorem~\ref{thm:main}.
For $i=0$, both sides of the asserted bound equal $1$.
For $1\le i\le d$, equality forces $B_1=H_1$, because all summands above are nonnegative and $\binom{d-1}{i-1}>0$.
The equality argument in the proof of Theorem~\ref{thm:main} then implies that $G$ is a cactus graph.
Conversely, for a cactus graph Theorem~\ref{thm:main} gives $B_k=H_k$ for every $k$, so equality holds for every $i$.

In particular, setting $i=d$ and using $F_d=\tau(G)$ yields $\tau(G)\le\prod\limits_{j=1}^{d}c_j$.
When $d\ge1$, apply the equality characterization for $F_i$ above with $i=d$. When $d=0$, $G$ is a tree, $\tau(G)=1$, and the empty product is $1$, so the conclusion also holds.
\end{proof}

\section{Conclusion}\label{sec:conclusion}

Using a quotient algebra representation associated with the cographic matroid, we establish upper bounds on the coefficients $H_i$ of the all-terminal reliability polynomial in terms of the lengths of the cycles in an arbitrary binary cycle basis, thereby providing an affirmative answer to the aforementioned question for binary cycle bases. Through the coefficient transformation relating $F_i$ and $H_i$, we further derive upper bounds on the numbers of deletable edge sets and spanning trees. We also prove that the upper bounds on $H_i$ are attained simultaneously for all $0\le i\le d$ if and only if $G$ is a cactus graph.
The bounds obtained in this paper depend only on the lengths of the cycles in a cycle basis. Several questions remain for further investigation. How can the number and distribution of edges shared by cycles in a cycle basis be used to obtain sharper upper bounds on the coefficients $H_i$? How can the lengths of the basis cycles and the cutset structure of the graph be used together to improve estimates for the numbers of deletable edge sets and spanning trees? When edges have different failure probabilities, how can the upper bounds on reliability obtained by the cycle-basis method be further sharpened?

\noindent{\bf Conflicts of Interest:}\\
 The authors declare that they have no conflicts of interest.

\noindent{\bf Data availability statement:} \\
No data were used to support this study.

\renewcommand\refname{\ References}

\end{document}